\documentclass[12pt,a4paper]{article}
\usepackage[a4paper,margin=15mm]{geometry}
\usepackage[T1]{fontenc}
\usepackage{amsmath,amssymb,amsthm}
\usepackage{microtype}
\usepackage[hidelinks]{hyperref}
\newtheorem{theorem}{Theorem}
\newtheorem{lemma}[theorem]{Lemma}
\newtheorem*{conjecture}{Conjecture}
\newtheorem{known}{Theorem}

\newcommand{\Prob}{\mathbb P}

\title{A proof of a conjecture on bipartite correspondence packing}
\author{Xin Zhang\\
School of Mathematics and Statistics, Xidian University\\
Xi'an 710071, China\\
\texttt{xzhang@xidian.edu.cn}}
\date{}

\begin{document}
\maketitle
\begin{abstract}
A correspondence cover assigns a set of possible colors to each vertex of a graph. For each edge, a matching between the two endpoint color sets specifies the forbidden color pairs. A correspondence coloring selects one color at each vertex so that no forbidden pair is selected. A full correspondence packing is a collection of such colorings that uses every color at each vertex exactly once.

We prove a conjecture of Cambie, Cames van Batenburg, Davies, and Kang [\emph{Random Structures \& Algorithms} \textbf{64} (2024), 62--93]. For every sufficiently large $D$, there is a positive integer $k(D)=(1+o(1))D/\log D$ such that every correspondence cover of any finite bipartite graph has a full correspondence packing if each color has at most $D$ forbidden partners and each vertex has $k(D)$ colors. More precisely, we prove this conclusion with
\[
 k(D)=\left\lceil\frac{D}{\log D-2\log\log D-6}\right\rceil
\]
for every integer $D\ge\lceil e^{20}\rceil$.
\end{abstract}


\section{Introduction}

All graphs in this paper are finite and simple, with no loops or multiple edges. For a graph $G$, write $V(G)$ for its vertex set and $E(G)$ for its edge set. Two vertices are neighbors if they form an edge. The degree $\deg_G(x)$ is the number of neighbors of $x$, and $\Delta(G)$ is the largest vertex degree, taken to be zero when $G$ has no vertices. An independent set is a set of vertices containing no adjacent pair. A matching is a set of edges no two of which share an endpoint; it is perfect if every vertex is an endpoint of one of its edges. A graph is bipartite if its vertices can be divided into two independent sets.

Every base graph considered below has at least one vertex. Let $F$ be a base graph. For each $v\in V(F)$, let $L(v)$ be a nonempty finite set of possible colors for $v$. Colors at different base vertices may have the same label. We treat them as separate copies: a color labeled $1$ at $u$ and one labeled $1$ at $v$ are represented by $(u,1)$ and $(v,1)$. With this convention, the sets $L(v)$ are pairwise disjoint. Below, a color means one of these copies.

For each edge $uv\in E(F)$, choose a matching $M_{uv}$ between $L(u)$ and $L(v)$. Each edge of $M_{uv}$ has one endpoint in each of these two sets. The matching may be empty and need not cover either set. Define a graph $H$ by
\[
 V(H)=\bigcup_{v\in V(F)}L(v),
 \qquad
 E(H)=\bigcup_{uv\in E(F)}M_{uv}.
\]
We call $H$, together with the specified sets $L(v)$, a correspondence cover of $F$. Each color is a vertex of $H$, and the color sets $L(v)$ partition $V(H)$. These sets are also called the parts of the cover. By construction, each $L(v)$ is an independent set in $H$, and there are no edges between $L(u)$ and $L(v)$ when $uv\notin E(F)$. Also, each color has at most one neighbor in any other color set, because the edges between two color sets form a matching.

For distinct base vertices $u$ and $v$, assigning $c\in L(u)$ to $u$
and $c'\in L(v)$ to $v$ creates a conflict exactly when $cc'\in E(H)$.
In particular, two colors with the same label do not conflict unless they
are joined by an edge of $H$. A correspondence coloring assigns one color from $L(v)$ to each base
vertex $v$, so that no two chosen colors are joined by an edge of $H$.
The chosen colors form an independent transversal, that is, an independent set
of $H$ containing exactly one color from each $L(v)$.

For a positive integer $k$, the cover is $k$-fold if $|L(v)|=k$ for every $v\in V(F)$, where $|S|$ denotes the size of a finite set $S$. A full correspondence packing is a partition of $V(H)$ into independent transversals. Each transversal contains one of the $k$ colors in any fixed $L(v)$. Thus, a full packing consists of $k$ colorings. For each vertex $v$, these colorings use every element of $L(v)$ exactly once.

The maximum degree $\Delta(F)$ counts the largest number of neighbors of a base vertex. The maximum degree $\Delta(H)$ counts the largest number of forbidden partners of a single color. If $V(F)=A\cup B$ is a bipartition of $F$, put
\[
 U_A=\bigcup_{v\in A}L(v),\qquad
 U_B=\bigcup_{v\in B}L(v).
\]
Every edge of $H$ belongs to a matching $M_{uv}$ for an edge $uv$ of $F$. Its endpoints therefore lie in different sets $U_A,U_B$. These sets form a bipartition of $H$.

All logarithms are natural.

Cambie, Cames van Batenburg, Davies, and Kang proved the following result with a bound on the maximum degree of the base graph.

\begin{samepage}
\begin{known}[Cambie, Cames van Batenburg, Davies, and Kang {\cite[Theorem 11]{CCDK}}]\label{thm:known}
For every $\varepsilon>0$, there is $D_0>1$ such that the following holds for every integer $D\ge D_0$. Let $F$ be a finite bipartite graph with
\[
 \Delta(F)\le D.
\]
Then every $k$-fold correspondence cover of $F$, where
\[
 k=\left\lceil(1+\varepsilon)\frac{D}{\log D}\right\rceil,
\]
has a full correspondence packing.
\end{known}
\end{samepage}

For a color $c\in L(v)$, every neighbor of $c$ in $H$ lies in $L(u)$ for some neighbor $u$ of $v$ in $F$. Each such color set $L(u)$ contains at most one neighbor of $c$ in $H$. Consequently,
\[
 \deg_H(c)\le\deg_F(v),\qquad \Delta(H)\le\Delta(F).
\]
Thus the condition $\Delta(F)\le D$ in Theorem~\ref{thm:known} implies $\Delta(H)\le D$. The converse need not hold, even when every edge of the base graph has a nonempty matching in the cover.

For example, fix integers $k\ge2$ and $D\ge1$. Take a base graph consisting of one vertex $v$ and $kD$ neighbors $v_{i,j}$, where $1\le i\le k$ and $1\le j\le D$, with no other edges. Give every base vertex a color set of size $k$, and write $L(v)=\{c_1,\ldots,c_k\}$. For each $i,j$, choose $c'_{i,j}\in L(v_{i,j})$ and set
\[
 M_{v v_{i,j}}=\{c_i c'_{i,j}\}.
\]
In the base graph, $v$ has degree $kD$ and every other vertex has degree one. In the cover, each $c_i$ has degree $D$, each $c'_{i,j}$ has degree one, and all remaining colors have degree zero in $H$. Hence
\[
 \Delta(F)=kD>D,\qquad \Delta(H)=D.
\]
The condition on the cover degree therefore includes covers whose base graphs fail the degree condition in Theorem~\ref{thm:known}.

Cambie, Cames van Batenburg, Davies, and Kang~\cite{CCDK} asked whether a bound on $\Delta(H)$ alone is sufficient at the scale $D/\log D$. In the following conjecture, $k(D)$ depends only on $D$; no bound is imposed on the degree or the number of vertices of the finite bipartite base graph.

\begin{conjecture}[Cambie, Cames van Batenburg, Davies, and Kang {\cite[Conjecture 14]{CCDK}}]
There are positive integers $k(D)$ satisfying
\[
 k(D)\le(1+o(1))\frac{D}{\log D}\qquad(D\longrightarrow\infty)
\]
such that every $k(D)$-fold correspondence cover $H$ of a finite bipartite graph $F$ with
\[
 \Delta(H)\le D
\]
has a full correspondence packing.
\end{conjecture}

We prove this conjecture with an explicit choice of $k(D)$.

\begin{samepage}
\begin{theorem}\label{thm:main}
For each integer $D\ge\lceil e^{20}\rceil$, set
\[
 k(D)=\left\lceil\frac{D}{\log D-2\log\log D-6}\right\rceil.
\]
Let $H$ be a correspondence cover of a finite bipartite graph $F$ such that $\Delta(H)\le D$ and $|L(v)|=k(D)$ for every $v\in V(F)$. Then $H$ has $k(D)$ pairwise disjoint independent transversals, which form a full correspondence packing. Moreover,
\[
 k(D)=(1+o(1))\frac{D}{\log D}\qquad(D\longrightarrow\infty).
\]
\end{theorem}
\end{samepage}

To define $k(D)$ for smaller nonnegative integers $D$, put $D_0=\lceil e^{20}\rceil$ and set $k(D)=k(D_0)$ when $D<D_0$. If $\Delta(H)\le D<D_0$, Theorem~\ref{thm:main} applied with $D_0$ gives a full packing for every $k(D_0)$-fold cover. These finitely many values do not affect the asymptotic formula for $k(D)$.

Our main result gives a sufficient condition for partitioning
the vertices of $H$ into independent transversals. It fixes the size $n$ of each color set $L(v)$ and gives a sufficient upper bound on the cover degree.

\begin{theorem}\label{thm:finite}
Let $n\ge\lceil e^{10}\rceil=22027$ be an integer. Let $H$ be a correspondence cover of a finite bipartite graph $F$, with $|L(v)|=n$ for all $v\in V(F)$. If
\begin{equation}\label{eq:main}
 \Delta(H)\le n(\log n-\log\log n-6),
\end{equation}
then $V(H)$ can be partitioned into $n$ independent transversals.
\end{theorem}

We prove Theorem~\ref{thm:finite} in Section~\ref{sec:finite-proof}.
In Section~\ref{sec:main-proof}, we apply it with $n=k(D)$ to prove
Theorem~\ref{thm:main}.

\subsection{Related results}

We follow the random ordering and matching framework of Cambie,
Cames van Batenburg, Davies, and Kang \cite[Section~7]{CCDK}.
Their proof completes the matchings between color sets and uses a
bound on $\Delta(F)$. Completing these matchings can increase
$\Delta(H)$, so we work directly with the partial matchings.
In the availability matrix defined in Section~\ref{sec:setup},
a color set can make an entry in the column for $c$ zero only if it
contains a neighbor of $c$ in $H$. Thus at most $\deg_H(c)$ color sets
can affect that entry. We use this observation to bound the probability
of an all-zero submatrix in terms of $\Delta(H)$ and obtain a full
packing.

For fixed $\varepsilon>0$ and sufficiently large $D$, Cambie and Kang
\cite[Corollary~1.3]{CK} proved that a correspondence cover $H$ of a
bipartite graph $F$ with $\Delta(H)\le D$ has an independent transversal
if every color set has at least $(1+\varepsilon)D/\log D$ colors. Cambie, Haxell, Kang, and
Wdowinski \cite{CHKW} gave a sharp asymmetric condition for an
independent transversal in a bipartite graph whose vertex classes each
lie entirely within one side of the bipartition.

Let $G$ be a graph whose vertices are partitioned into classes. Alon
\cite{Alon} proved that this partition has an independent transversal
if every class has at least $2e\Delta(G)$ vertices. Haxell
\cite{Haxell} improved this bound to $2\Delta(G)$.
The local degree of the partition is the largest number of neighbors
that a vertex has in any one class. Loh and Sudakov \cite{LS} showed
that, when the local degree is $o(\Delta(G))$, classes of size at least
$\Delta(G)+o(\Delta(G))$ guarantee an independent transversal. Glock
and Sudakov \cite{GS}, and Kang and Kelly \cite{KK}, obtained bounds
using the average degrees of vertices in the classes.

For packing, Loh and Sudakov \cite{LS} conjectured that all vertices can be partitioned into independent transversals under their local degree assumptions. Chakraborti and Tran \cite{CT} proved an approximate form: for fixed $\varepsilon\in(0,1)$ and sufficiently large $n$, classes of size $n$, maximum degree at most $(1-\varepsilon)n$, and sufficiently small local degree relative to $n$ guarantee at least $\lfloor(1-\varepsilon)n\rfloor$ disjoint transversals. Yuster \cite{Yuster} proved a full partition when there are $r$ classes of size $n\ge1.78r$, the edges between each pair of classes form a matching, and $r$ is sufficiently large.

List coloring is the special case in which colors at different vertices have shared labels and the forbidden pairs across each edge are precisely the pairs with equal labels. A list packing uses every element of every vertex's list once in a collection of proper list colorings. Further results treat bounded base degree \cite{Bounded}, complete bipartite graphs \cite{Complete}, and planar graphs \cite{Planar}.  Fractional packing asks for a probability distribution on colorings in which each color in $L(v)$ is chosen with probability $1/|L(v)|$ \cite{Fractional}. Approximate counting and sampling of list packings are studied in \cite{Sampling}.

\subsection{Proof outline}

We outline the proof of Theorem~\ref{thm:finite}. As in \cite{CCDK},
we use random orderings and auxiliary matchings.
Let $A,B$ be the two sides of the base graph $F$, and suppose that each color set $L(v)$ has size $n$. Write $[n]=\{1,\ldots,n\}$. We construct $n$ colorings, represented by rows $i\in[n]$. For each $a\in A$, choose a uniformly random ordering of the colors in $L(a)$, independently for different base vertices $a$, and place its $i$-th color in row $i$. Each row now assigns one color to every vertex of $A$. These assignments create no conflicts as there are no edges of $F$ within $A$, so there are no edges of $H$ between the color sets of distinct vertices of $A$. We will complete each row by assigning colors to the vertices of $B$.

Fix $b\in B$ and put $P=L(b)$. Form an auxiliary bipartite graph $J_P$ with sides $[n]$ and $P$. Join row $i$ to color $c\in P$ when $c$ has no neighbor in $H$ among the colors already in that row. Such a color is available in row $i$, meaning that it can be assigned to $b$ in that row without a conflict. A perfect matching of $J_P$ is therefore equivalent to assigning all colors in $P$ to distinct allowed rows. If every $J_P$ has a perfect matching, we can make these assignments for every $b\in B$. There are no edges of $H$ between the color sets of distinct vertices of $B$, so the completed rows are valid colorings. Each color is used in exactly one row, giving a full packing. Equivalence~\eqref{eq:packing-matchings} in Section~\ref{sec:setup} states this relation for the fixed orderings on the $A$ side.

For this fixed $b$, let $M_P$ be the $n\times n$ matrix with rows indexed by $i\in[n]$ and columns by colors $c\in P=L(b)$. Set $M_P(i,c)=1$ if $c$ is available in row $i$, and $M_P(i,c)=0$ otherwise. Thus the matrix is associated with the base vertex $b$, and each column corresponds to one color in $L(b)$. For $S\subseteq[n]$ and $X\subseteq P$, write $M_P[S,X]$ for the submatrix with rows in $S$ and columns corresponding to colors in $X$. Lemma~\ref{lem:hall-obstruction}, proved using Hall's theorem, gives
\[
 \begin{aligned}
 &J_P\text{ has no perfect matching}\\
 &\quad\Longleftrightarrow\quad
 \begin{gathered}
 \text{there exist }S\subseteq[n]\text{ and }X\subseteq P\text{ with }|S|+|X|=n+1\\
 \text{such that every entry of }M_P[S,X]\text{ is zero}.
 \end{gathered}
 \end{aligned}
\]
This equivalence reduces matching failure to the occurrence of one of these all-zero submatrices.

Fix such sets $S$ and $X$. We consider selected columns one at a time when $|X|\le|S|$, and selected rows one at a time otherwise. If $M_P[S,X]$ is all zero, then these selected columns or rows must all be zero. 
We can therefore bound the probability that the whole submatrix is all
zero by considering at most $h=\lfloor n/\log n\rfloor$ selected
columns or rows, one at a time. We next explain how conflicts make a
column or a row of $M_P[S,X]$ all zero.

For a color $c\in X$, its column records the rows in which $c$ can be assigned to $b$. All entries of this column in the rows $S$ are zero exactly when every row in $S$ already contains a color adjacent to $c$ in $H$. To determine this event, we consider the row containing each neighbor of $c$. For each base vertex $a\in A$, consider its color set $L(a)$. The edges of $H$ between $L(a)$ and $P=L(b)$ form a matching, so $c$ has at most one neighbor in $L(a)$. If this neighbor exists, it lies in exactly one row and makes $c$ unavailable there. If there is no such neighbor, $L(a)$ creates no conflict with $c$. Thus $L(a)$ can make at most one entry in this column zero. In the occupancy model of Lemma~\ref{lem:occupancy}, the bins are the rows of $S$. For each color set $L(a)$ containing a neighbor of $c$, the random row of that neighbor defines one trial: it selects the bin for that row if the row is in $S$, and selects no bin otherwise. The column is all zero on $S$ exactly when every bin is occupied.

For a selected row $i$, we instead consider the color in row $i$ of each color set $L(a)$, where $a\in A$. All entries in the columns $X$ are zero exactly when every color in $X$ has a neighbor among these chosen colors. Each chosen color has at most one neighbor in $P$, again by the matching condition. Now the bins are the colors of $X$, and the trial from $L(a)$ selects the neighbor in $X$ of its chosen color, if there is one. Thus every bin is occupied exactly when the row is all zero on $X$.

At step $j$, where $1\le j\le h$, we consider the $j$th selected
column or row. When
considering a column, we condition on the rows containing the
neighbors in $H$ of the colors indexing the first $j-1$ selected
columns. When considering a row, we condition on the color from
each color set over $A$ in each of the first $j-1$ selected rows.
We consider every such choice with positive probability.
By Lemma~\ref{lem:conditional-orderings}, conditional on these
choices, the remaining colors in each color set are ordered
uniformly among the remaining rows. These orderings are independent
across the color sets.
Put $D=\Delta(H)$. At most $D$ trials can select any fixed bin, and each does so with conditional probability at most $1/(n-h+1)$. By independence, the conditional probability that the bin remains unoccupied is at least $(1-1/(n-h+1))^D$. Lemma~\ref{lem:occupancy} converts these bounds into an upper bound on the probability that the next selected column or row is all zero. Lemma~\ref{lem:conditional-steps} combines these bounds by summing over the possible earlier values. This gives the bound for an all-zero submatrix in Lemma~\ref{lem:rectangle}.

Let $E_P$ be the event that $J_P$ has no perfect matching. Using the counts in Lemma~\ref{lem:submatrix-count}, we sum the submatrix bounds over pairs $(S,X)$ with $|S|+|X|=n+1$ to bound $\Prob(E_P)$. The event $E_P$ depends only on the orderings of color sets $L(a)$ joined to $P$ by an edge of $H$. Join two events whenever they use a common ordering. Each color set has $n$ colors of degree at most $D$, so it is joined to at most $nD$ other color sets. There are therefore at most $(nD)^2$ events joined to any one event. Lemma~\ref{lem:matching-dependencies} verifies that this is a dependency graph, including the required independence from the joint collection of nonneighbor events.

Under the hypothesis of Theorem~\ref{thm:finite}, the probability
bound and the dependency bound satisfy the condition of the local
lemma. Hence, with positive probability, none of the events $E_P$
occurs. Every auxiliary graph then has a perfect matching, so the
rows can be completed to a full packing. This proves
Theorem~\ref{thm:finite}. 

\section{Random orderings and matchings}\label{sec:setup}

We express the required color assignments as perfect matchings in auxiliary graphs. We then describe matching failure by all-zero submatrices and count the possible pairs of row and color sets. Write $\Prob$ for probability.

Let $H$ be an $n$-fold correspondence cover of a bipartite base graph $F$. Choose a bipartition $V(F)=A\cup B$, with no edges inside either side. Write
\[
 \mathcal A=\{L(a):a\in A\},\qquad
 \mathcal B=\{L(b):b\in B\}
\]
for the two collections of color sets. Each member of $\mathcal A$ or $\mathcal B$ is one set $L(v)$, whereas each member of $A$ or $B$ is a base vertex. There are no edges of $H$ between color sets in the same collection.

Write $[n]=\{1,\ldots,n\}$. For each $Q\in\mathcal A$, choose independently a uniform bijection $\pi_Q:[n]\to Q$; that is, each of its $n!$ orderings has the same probability. The color from $Q$ placed in row $i$ is $\pi_Q(i)$. Thus row $i$ currently contains $\{\pi_Q(i):Q\in\mathcal A\}$. This is an independent set in $H$, since there are no edges between color sets in $\mathcal A$. We seek to add one color from each color set in $\mathcal B$ to this row.

Fix a color set $P=L(b)$ with $b\in B$, so $P\in\mathcal B$. A color $c\in P$ is available in row $i$ if it has no neighbor in $H$ among the colors already in that row. Define the availability graph $J_P$ to be the bipartite graph with sides $[n]$ and $P$, putting an edge between $i$ and $c$ precisely when
\[
 \pi_Q(i)c\notin E(H)\quad\text{for every }Q\in\mathcal A.
\]
The condition ranges over every $Q$ because $c$ must be nonadjacent in $H$ to every color already in row $i$. It holds automatically for any $Q$ with no edge of $H$ between $Q$ and $P$. An edge of $J_P$ represents an allowed assignment; an edge of $H$ represents a conflict between two colors.

Fix the orderings $\pi_Q$ for all $Q\in\mathcal A$. For these fixed orderings,
\begin{equation}\label{eq:packing-matchings}
 \begin{aligned}
 &\text{the rows can be completed to a full packing}\\
 &\qquad\Longleftrightarrow\quad
 J_P\text{ has a perfect matching for every }P\in\mathcal B.
 \end{aligned}
\end{equation}
Indeed, a perfect matching of $J_P$ assigns each color of $P$ to a distinct row in which it has no neighbor among the colors already chosen on the $A$ side. There are no edges of $H$ between color sets over the $B$ side, so these assignments can be made for all $P\in\mathcal B$ without creating any further conflicts. Each completed row contains one color from every color set, and every color is used in exactly one row. Thus the completed rows form a full packing. Conversely, a full packing that completes the fixed rows assigns every color of each $P\in\mathcal B$ to exactly one row. The colors in each completed row form an independent set in $H$,
so each assignment corresponds to an edge of $J_P$. These edges use each row and each color of $P$ exactly once, so they form a perfect matching of $J_P$.

Define the availability matrix $M_P$, with rows indexed by $[n]$ and columns indexed by $P$, by
\[
 M_P(i,c)=
 \begin{cases}
 1,&\text{if }ic\in E(J_P),\\
 0,&\text{otherwise}.
 \end{cases}
\]
Thus an entry is one precisely when the corresponding color is available in the corresponding row. For nonempty sets $S\subseteq[n]$ and $X\subseteq P$, write $M_P[S,X]$ for the submatrix obtained by selecting the rows in $S$ and the columns in $X$. Define $\mathcal R(S,X)$ to be the event that every entry of $M_P[S,X]$ is zero. Equivalently, every color in $X$ is unavailable in every row of $S$. By the definition of $J_P$, this means that for every $i\in S$ and $c\in X$ there is some $Q\in\mathcal A$ with $\pi_Q(i)c\in E(H)$. The color set $Q$ can be different for different pairs $(i,c)$.

We use Hall's theorem in the following form.

\begin{samepage}
\begin{lemma}[Hall's theorem; see {\cite[Theorem 16]{CCDK}}]\label{lem:hall}
Let $J$ be a bipartite graph with sides $U,W$ of equal finite size. For $S\subseteq U$, let $N_J(S)$ be the set of vertices of $W$ adjacent to at least one vertex of $S$. Then $J$ has a perfect matching if and only if
\[
 |N_J(S)|\ge |S|\qquad\text{for every }S\subseteq U.
\]
\end{lemma}
\end{samepage}

\begin{lemma}\label{lem:hall-obstruction}
For each $P\in\mathcal B$, the graph $J_P$ has no perfect matching if and only if there are nonempty sets $S\subseteq[n]$ and $X\subseteq P$ such that $|S|+|X|=n+1$ and every entry of $M_P[S,X]$ is zero. Equivalently,
\begin{equation}\label{eq:hall-failure-events}
 \{J_P\text{ has no perfect matching}\}
 =\bigcup_{\substack{S\subseteq[n],\ X\subseteq P\\|S|+|X|=n+1}}\mathcal R(S,X).
\end{equation}
\end{lemma}

\begin{proof}
Apply Lemma~\ref{lem:hall} to $J_P$, whose sides are $[n]$ and $P$. For $S\subseteq[n]$, the set $N_{J_P}(S)$ consists of the colors in $P$ available in at least one row of $S$. If $J_P$ has no perfect matching, the lemma gives a nonempty $S\subseteq[n]$ with
\[
 |N_{J_P}(S)|\le |S|-1.
\]
Writing $s=|S|$, we have
\[
 |P\setminus N_{J_P}(S)|=n-|N_{J_P}(S)|\ge n-s+1.
\]
Choose $X\subseteq P\setminus N_{J_P}(S)$ with $x=|X|=n-s+1$. For every $c\in X$, being outside $N_{J_P}(S)$ means that $ic\notin E(J_P)$ for every $i\in S$. Thus every entry of $M_P[S,X]$ is zero, so $\mathcal R(S,X)$ holds, with
\begin{equation}\label{eq:hall-sizes}
 s+x=n+1.
\end{equation}

Conversely, suppose that $M_P[S,X]$ is all zero and $s+x=n+1$. No color in $X$ is adjacent in $J_P$ to a row in $S$, so $N_{J_P}(S)\subseteq P\setminus X$. Consequently,
\[
 |N_{J_P}(S)|\le n-x=s-1<s.
\]
Hall's condition in Lemma~\ref{lem:hall} fails, so $J_P$ has no perfect matching. This proves both directions. Since the equivalence holds for every choice of the orderings, it also gives the event equality \eqref{eq:hall-failure-events}.
\end{proof}

We next count the pairs of sets in Lemma~\ref{lem:hall-obstruction}. The count concerns all possible pairs, regardless of whether their submatrices are zero.

\begin{lemma}\label{lem:submatrix-count}
Let $P$ be a set of size $n$. For each integer $\mu$ with $1\le \mu\le\lfloor(n+1)/2\rfloor$, there are at most $2n^{2\mu-1}$ pairs $S\subseteq[n]$, $X\subseteq P$ satisfying
\[
 |S|+|X|=n+1,\qquad \min\{|S|,|X|\}=\mu.
\]
There are at most $4^n$ pairs satisfying $|S|+|X|=n+1$ in total.
\end{lemma}

\begin{proof}
Write $s=|S|$ and $x=|X|$, so $s+x=n+1$.

For fixed sizes $s,x$, there are $\binom ns$ choices of $S$ and, for each such choice, $\binom nx$ choices of $X$.
Thus there are $\binom ns\binom nx$ pairs $(S,X)$. Put $\mu=\min\{s,x\}$. If $s\le x$, then $s=\mu$ and $x=n+1-\mu$. If $x<s$, then $x=\mu$ and $s=n+1-\mu$. In either case, we have
\begin{equation}\label{eq:count}
 \binom ns\binom nx
 =\binom n\mu\binom n{n+1-\mu}
 =\binom n\mu\binom n{\mu-1}
 \le n^\mu n^{\mu-1}=n^{2\mu-1}.
\end{equation}
The second equality follows by taking complements, and the inequality
uses $\binom nr\le n^r$. Since $\mu=s$ or $\mu=x$, there are at most two possible ordered size
pairs. Thus there are at most $2n^{2\mu-1}$ pairs $(S,X)$ for each $\mu$.

For the total count, each element can either belong to a subset or not, so each side has $2^n$ subsets. Hence there are $2^n\cdot2^n=4^n$ pairs of subsets altogether. This count includes pairs that fail $s+x=n+1$, and therefore bounds the required number of pairs as well.
\end{proof}

\section{Probability estimates}\label{sec:occupancy}

The next lemma describes the orderings after some colors have been
assigned to specified rows.

\begin{lemma}\label{lem:conditional-orderings}
Use the independent uniform orderings $\pi_Q:[n]\to Q$ from
Section~\ref{sec:setup}. For each $Q\in\mathcal A$, assign $r_Q$
distinct colors of $Q$ to distinct rows, where $0\le r_Q\le n$. Let
$K$ be the event that all these assignments hold. Then $\Prob(K)>0$.
Conditional on $K$, the orderings are independent across the color
sets. In each $Q$, the remaining colors are ordered uniformly among
the remaining rows. In particular, if $r_Q<n$, then for every
remaining row $i$ and every unassigned color $c\in Q$,
\[
 \Prob(\pi_Q(i)=c\mid K)=\frac1{n-r_Q}.
\]
\end{lemma}

\begin{proof}
For each $Q\in\mathcal A$, let $K_Q$ be the event that each of the
$r_Q$ specified colors in $Q$ is placed in its assigned row by
$\pi_Q$. Thus $K=\bigcap_{Q\in\mathcal A}K_Q$.
Fix $Q\in\mathcal A$. For $K_Q$ to occur, the $r_Q$ specified colors
must occupy their assigned rows. The remaining $n-r_Q$ colors can
then be placed in the remaining $n-r_Q$ rows in exactly
$(n-r_Q)!$ ways. Hence
\[
 \Prob(K_Q)=\frac{(n-r_Q)!}{n!}>0.
\]
Since the orderings are independent and $K=\bigcap_{Q\in\mathcal A}K_Q$, we have $\Prob(K)=\prod_{Q\in\mathcal A}\Prob(K_Q)>0$.

Choose one ordering $\rho_Q$ satisfying $K_Q$ for each $Q$. Independence before conditioning gives
\[
 \begin{aligned}
 \Prob(\pi_Q=\rho_Q\text{ for all }Q\in\mathcal A\mid K)
 &=\frac{\prod_{Q\in\mathcal A}(1/n!)}{\prod_{Q\in\mathcal A}((n-r_Q)!/n!)}\\
 &=\prod_{Q\in\mathcal A}\frac1{(n-r_Q)!}.
 \end{aligned}
\]
Thus the conditional distribution is the product of uniform distributions on the permitted orderings. This proves the claimed conditional independence and uniformity. If $r_Q<n$, a specified unassigned color occupies a specified remaining row in exactly $(n-r_Q-1)!$ of the $(n-r_Q)!$ permitted orderings. The ratio of these counts is $1/(n-r_Q)$, proving the remaining assertions.
\end{proof}

The next lemma combines probability bounds from successive steps.
In our applications of this lemma, the events $K$ represent possible
choices of the values fixed in earlier steps.

\begin{lemma}\label{lem:conditional-steps}
Let $C_1,\ldots,C_k$ be events. Put $B_0=\Omega$, where $\Omega$ is
the sample space, and
\[
 B_j=C_1\cap\cdots\cap C_j
 \qquad (1\le j\le k).
\]

For each $j$, let $\mathcal K_j$ be a finite collection of disjoint
events of positive probability whose union is $\Omega$. Suppose
that each $K\in\mathcal K_j$ either lies in $B_{j-1}$ or is
disjoint from $B_{j-1}$.

Let $a_j\in[0,1]$. If
\[
 \Prob(C_j\mid K)\le a_j
\]
for every $K\in\mathcal K_j$ contained in $B_{j-1}$, then
\[
 \Prob(B_j)\le a_j\Prob(B_{j-1})
 \qquad (1\le j\le k).
\]
Consequently,
\[
 \Prob(C_1\cap\cdots\cap C_k)\le\prod_{j=1}^k a_j.
\]
\end{lemma}

\begin{proof}
Fix $j$. The events in $\mathcal K_j$ partition $\Omega$, and each
is either contained in $B_{j-1}$ or disjoint from it. Hence those
contained in $B_{j-1}$ partition $B_{j-1}$. For each such $K$,
since $B_j=B_{j-1}\cap C_j$, we have
\[
 B_j\cap K=C_j\cap K.
\]

In the sums below, $K$ ranges over the events in $\mathcal K_j$
contained in $B_{j-1}$. Since these events partition $B_{j-1}$ and
$B_j\subseteq B_{j-1}$, the events $B_j\cap K$ are disjoint and
their union is $B_j$. 
Therefore
\[
 \begin{aligned}
 \Prob(B_j)
 &=\sum_{K\subseteq B_{j-1}}\Prob(B_j\cap K)=\sum_{K\subseteq B_{j-1}}\Prob(C_j\cap K)
  =\sum_{K\subseteq B_{j-1}}\Prob(K)\Prob(C_j\mid K)\\
 &\le a_j\sum_{K\subseteq B_{j-1}}\Prob(K)
  =a_j\Prob(B_{j-1}).
 \end{aligned}
\]

Starting from $\Prob(B_0)=1$ and applying this inequality for
$j=1,\ldots,k$ gives
\[
 \Prob(B_k)\le\prod_{j=1}^k a_j.
\]
Since $B_k=C_1\cap\cdots\cap C_k$, this proves the lemma.
\end{proof}

Dubhashi and Ranjan \cite[Theorem~13]{DR} showed that the bin counts
in their balls-and-bins model are negatively associated. We use the
following consequence to bound the conditional probabilities in
Lemma~\ref{lem:conditional-steps}.

\begin{lemma}[Dubhashi and Ranjan]\label{lem:nonempty-bins}
Let a finite collection of balls be placed independently into a
nonempty finite set $J$ of bins. Each ball is placed in exactly one
bin, and the placement distributions may differ between balls.
For $i\in J$, let $B_i$ be the number of balls in bin $i$. Then, for
every $I\subseteq J$,
\[
 \Prob(B_i\ge1\text{ for every }i\in I)
 \le \prod_{i\in I}\Prob(B_i\ge1).
\]
\end{lemma}

In an occupancy experiment, the elements of a finite set $I$ are
called bins. Each trial selects at most one bin and may select none.
A bin is occupied if at least one trial selects it.

\begin{lemma}\label{lem:occupancy}
Let a finite collection of independent trials each select at most
one bin in $I$. For $i\in I$, let $q_i$ be the probability that no
trial selects $i$. Then
\[
 \Prob(\text{every bin in }I\text{ is occupied})
 \le \prod_{i\in I}(1-q_i)
 \le \exp\left(-\sum_{i\in I}q_i\right).
\]
\end{lemma}

\begin{proof}
Add an auxiliary bin $\ast\notin I$, and put $J=I\cup\{\ast\}$.
For each trial, place a ball in the selected bin, or in $\ast$ if
the trial selects no bin in $I$. These placements are independent,
and each ball is placed in exactly one bin of $J$.

For $i\in I$, let $B_i$ be the number of balls in bin $i$. This bin
is empty exactly when no trial selects it, so
$\Prob(B_i=0)=q_i$. Applying Lemma~\ref{lem:nonempty-bins} to the
set of bins $J$ and its subset $I$ gives
\[
 \begin{aligned}
 \Prob(\text{every bin in }I\text{ is occupied})
 &=\Prob(B_i\ge1\text{ for every }i\in I)\\
 &\le\prod_{i\in I}\Prob(B_i\ge1)
 =\prod_{i\in I}(1-q_i).
 \end{aligned}
\]
The second inequality follows from $1-q_i\le e^{-q_i}$.
\end{proof}

\section{All-zero submatrices}\label{sec:rectangles}

Use the random orderings from Section~\ref{sec:setup}. Recall that $A,B$ are the two sides of the base graph $F$, and $\mathcal B=\{L(b):b\in B\}$ is the collection of color sets over $B$. Fix $P=L(b)\in\mathcal B$. The matrix $M_P$ has one row for each $i\in[n]$ and one column for each color $c\in P$. Its entry $M_P(i,c)$ is one if $c$ is available in row $i$, and zero otherwise.

Choose nonempty sets $S\subseteq[n]$ and $X\subseteq P$, and put $s=|S|$ and $x=|X|$. Recall that $\mathcal R(S,X)$ is the event that all entries of the submatrix $M_P[S,X]$ are zero. The next lemma bounds its probability for every such pair $S,X$.

\begin{lemma}\label{lem:rectangle}
Suppose that every color in $P$ has degree at most $D$ in $H$. Let $h$ be an integer with $1\le h\le\lfloor n/2\rfloor$, and set
\[
 \tau=\left(1-\frac1{n-h+1}\right)^D.
\]
For $\mu=\min\{s,x\}$ and $\nu=\max\{s,x\}$, we have
\begin{equation}\label{eq:rectangle}
 \Prob(\mathcal R(S,X))\le\exp\bigl(-\nu\min\{\mu,h\}\tau\bigr).
\end{equation}
\end{lemma}

\begin{proof}
Put $k=\min\{\mu,h\}$. When $x\le s$, we select $k$ columns of
$M_P[S,X]$; when $s<x$, we select $k$ rows. We consider the
selected columns or rows one at a time. In either case, each selected
column or row has $\nu=\max\{s,x\}$ entries in the submatrix.

Recall that $\mathcal A=\{L(a):a\in A\}$ is the collection of color
sets over the $A$ side of the base graph. For $Q\in\mathcal A$,
the color $\pi_Q(i)$ is the color of $Q$ in row $i$, and
$\pi_Q^{-1}(c')$ is the row containing $c'\in Q$.

For $1\le j\le k$, let $C_j$ be the event that every entry of the
$j$th selected column or row in $M_P[S,X]$ is zero. We first describe
the two choices of columns or rows and the corresponding trials.

\smallskip
\noindent\textit{Case 1: $x\le s$.}

Choose distinct colors $c_1,\ldots,c_k\in X$, where
$k=\min\{x,h\}$, independently of the random orderings. We consider
the corresponding columns in this order. Thus $C_j$ is the event
that $c_j$ is unavailable in every row of $S$.

Fix $j$. For each color set $Q=L(a)$ with $a\in A$, consider the
colors in $Q$ adjacent in $H$ to at least one of
$c_1,\ldots,c_{j-1}$. Assign a row to each of these colors, using
distinct rows for distinct colors in the same set $Q$. Colors in
different sets may be assigned to the same row. Choose any such
assignment in every color set, and let $K$ be the event that all
these assignments hold in the random orderings.

Each assignment can be completed to an ordering of its color set,
so $\Prob(K)>0$. We calculate probabilities conditional on $K$.
When $j=1$, there are no assignments, and $K$ is the whole sample
space.

For each $Q\in\mathcal A$, let $r_Q$ be the number of colors of $Q$
assigned rows by $K$. By Lemma~\ref{lem:conditional-orderings},
conditional on $K$, the orderings in different color sets are
independent. In each $Q$, every color not assigned a row by $K$
is equally likely to occupy any of the remaining $n-r_Q$ rows.

The edges between $Q$ and $P$ form a matching. Hence $c_j$ has at
most one neighbor in $Q$. If this neighbor exists, it is distinct
from every neighbor in $Q$ of $c_1,\ldots,c_{j-1}$; otherwise it
would be incident with two edges of the matching. Thus no row was
assigned to this neighbor in the definition of $K$.

Each of $c_1,\ldots,c_{j-1}$ has at most one neighbor in $Q$, so
$r_Q\le j-1$. Since $j\le k\le h$,
\[
 n-r_Q\ge n-j+1\ge n-h+1.
\]

Use the rows in $S$ as bins. For each color set $Q\in\mathcal A$
containing a neighbor $c'$ of $c_j$, the random row
$\pi_Q^{-1}(c')$ defines one trial. The trial selects bin $i$ if
$c'$ lies in a row $i\in S$, and selects no bin if $c'$ lies
outside $S$.

Each trial depends only on its own ordering $\pi_Q$. Therefore,
conditional on $K$, the trials are independent, and each trial
selects at most one bin. There are exactly
$\deg_H(c_j)\le D$ trials: every neighbor of $c_j$ lies in a color
set in $\mathcal A$, and each such set contains at most one
neighbor.

Fix a row $i\in S$. For a color set $Q$ giving a trial, if $K$
assigns another color of $Q$ to row $i$, the trial selects $i$
with probability zero. Otherwise, the neighbor $c'$ of $c_j$ is
equally likely to occupy any of the $n-r_Q$ remaining rows.
Consequently,
\[
 \Prob(\text{the trial from $Q$ selects $i$}\mid K)
 \le \frac1{n-r_Q}\le\frac1{n-h+1}.
\]
The bin $i$ is empty exactly when no neighbor of $c_j$ lies in
row $i$. This is exactly the condition that $c_j$ is available
in row $i$. Thus $C_j$ occurs exactly when every bin in $S$ is
occupied.

\smallskip
\noindent\textit{Case 2: $s<x$.}

Choose distinct rows $i_1,\ldots,i_k\in S$, where
$k=\min\{s,h\}$, independently of the random orderings. We consider
these rows in this order. Thus $C_j$ is the event that every color
in $X$ is unavailable in row $i_j$.

Fix $j$. For each color set $Q=L(a)$ with $a\in A$, assign
$j-1$ distinct colors of $Q$ to the rows $i_1,\ldots,i_{j-1}$,
one color to each row. Choose any such assignment in every color
set, and let $K$ be the event that all these assignments hold in
the random orderings.

Each assignment can be completed to an ordering of its color set,
so $\Prob(K)>0$. We calculate probabilities conditional on $K$.
When $j=1$, there are no assignments, and $K$ is the whole sample
space.

In each $Q$, the event $K$ assigns $j-1$ colors to the distinct
rows $i_1,\ldots,i_{j-1}$. These colors cannot appear in row $i_j$.
By Lemma~\ref{lem:conditional-orderings}, conditional on $K$, the
orderings in different color sets are independent, and
$\pi_Q(i_j)$ is uniform among the remaining $n-j+1$ colors of $Q$.
Since $j\le k\le h$,
\[
 n-j+1\ge n-h+1.
\]

Use the colors in $X$ as bins. For each color set
$Q\in\mathcal A$, let $c'_Q=\pi_Q(i_j)$ be its color in row $i_j$.
This color defines one trial: if $c'_Q$ has a neighbor $c\in X$,
the trial selects bin $c$; otherwise it selects no bin.

The edges between $Q$ and $P$ form a matching, so $c'_Q$ has at
most one neighbor in $X\subseteq P$. Each trial depends only on
its own ordering $\pi_Q$. Therefore, conditional on $K$, the
trials are independent, and each trial selects at most one bin.

Fix a color $c\in X$. Every neighbor of $c$ lies in a color set
in $\mathcal A$, and each such set contains at most one neighbor
of $c$. Thus exactly $\deg_H(c)\le D$ color sets contain a
neighbor of $c$. Only the trials from these color sets can select
bin $c$.

For one of these sets $Q$, let $c'$ be the neighbor of $c$ in $Q$.
The trial from $Q$ selects $c$ exactly when $\pi_Q(i_j)=c'$. If
$K$ assigns $c'$ to one of the rows $i_1,\ldots,i_{j-1}$, then
the trial selects $c$ with probability zero. If $c'$ is not
assigned to an earlier row, it is one of the $n-j+1$ remaining
colors. Hence
\[
 \Prob(\text{the trial from $Q$ selects $c$}\mid K)
 \le \frac1{n-j+1}\le\frac1{n-h+1}.
\]
A color set containing no neighbor of $c$ selects $c$ with
probability zero.

The bin $c$ is empty exactly when none of the colors
$\pi_Q(i_j)$, with $Q\in\mathcal A$, is adjacent to $c$ in $H$.
This is exactly the condition that $c$ is available in row $i_j$.
Thus $C_j$ occurs exactly when every bin in $X$ is occupied.

\smallskip
\noindent\textit{The common probability bound.}

In either case, there are $\nu$ bins. Conditional on $K$, the trials
are independent and each trial selects at most one bin. Moreover,
for every bin, at most $D$ trials can select it, and each of these
trials selects it with conditional probability at most
$1/(n-h+1)$.

For a bin $z$, let $q_z$ be the conditional probability, given
$K$, that no trial selects $z$. By independence,
\[
 q_z\ge
 \left(1-\frac1{n-h+1}\right)^D
 =\tau.
\]
The event $C_j$ occurs exactly when all $\nu$ bins are occupied.
Applying Lemma~\ref{lem:occupancy} under the condition $K$ gives
\begin{equation}\label{eq:common-step-bound}
 \Prob(C_j\mid K)
 \le\exp\left(-\sum_z q_z\right)
 \le e^{-\nu\tau},
\end{equation}
where the sum is over the $\nu$ bins. This bound holds for every
choice of $K$ described in the relevant case.

Put $B_j=C_1\cap\cdots\cap C_j$, and let $B_0$ be the whole sample
space. For a fixed $j$, let $K$ range over the events obtained
from all assignments described in the relevant case. These events
partition the sample space.

Each such $K$ determines whether $B_{j-1}$ occurs. Indeed, in
Case 1, $K$ specifies the rows containing the neighbors of
$c_1,\ldots,c_{j-1}$; in Case 2, it specifies the color from every
color set in each of the rows $i_1,\ldots,i_{j-1}$. Thus each $K$
is contained in $B_{j-1}$ or disjoint from it.

These facts verify the partition assumption of
Lemma~\ref{lem:conditional-steps}. The bound
\eqref{eq:common-step-bound} verifies its conditional probability
assumption with $a_j=e^{-\nu\tau}$. The lemma gives
\[
 \Prob(B_j)\le e^{-\nu\tau}\Prob(B_{j-1}),
 \qquad
 \Prob(B_k)\le e^{-\nu k\tau}.
\]

If $\mathcal R(S,X)$ occurs, then every selected column or row is
all zero. Hence $\mathcal R(S,X)\subseteq B_k$, and therefore
\[
 \Prob(\mathcal R(S,X))
 \le\Prob(B_k)
 \le e^{-\nu k\tau}
 =\exp\bigl(-\nu\min\{\mu,h\}\tau\bigr).
\]
This proves \eqref{eq:rectangle}.
\end{proof}

\section{Proofs of the main theorems}\label{sec:finite}

We apply the estimate from Section~\ref{sec:rectangles} to prove Theorem~\ref{thm:finite}, then choose the color set size to deduce Theorem~\ref{thm:main}.

\subsection{The local lemma and dependencies}\label{sec:dependencies}

A dependency graph for a family of events has the events as its vertices. It has the property that each event is independent of the joint collection of events at all its nonneighbors, excluding itself. We use the following symmetric form of the Lov\'asz local lemma \cite[Section~1.1]{CK}.

\begin{samepage}
\begin{lemma}[Symmetric Lov\'asz local lemma]\label{lem:local}
Let $\mathcal E$ be a finite family of events, each of probability at most $p$, where $0\le p\le1$. Suppose that $\mathcal E$ has a dependency graph of maximum degree at most a nonnegative integer $d$. If
\[
 ep(d+1)\le1,
\]
then the probability that none of the events in $\mathcal E$ occurs is positive.
\end{lemma}
\end{samepage}

We next verify the dependency condition for the matching failures defined by our random orderings.

\begin{lemma}\label{lem:matching-dependencies}
Let $F$ be a bipartite graph with sides $A$ and $B$, and let $H$
be an $n$-fold correspondence cover of $F$ with $\Delta(H)\le D$.
Put
\[
 \mathcal A=\{L(a):a\in A\},
 \qquad
 \mathcal B=\{L(b):b\in B\}.
\]
For each $Q\in\mathcal A$, independently choose a uniform ordering
$\pi_Q:[n]\to Q$.

For $P\in\mathcal B$, let $J_P$ be the bipartite graph with sides
$[n]$ and $P$, where $i$ is joined to $c\in P$ if $c$ has no
neighbor in $H$ among the colors $\pi_Q(i)$, $Q\in\mathcal A$.
Let $E_P$ be the event that $J_P$ has no perfect matching.

Form a graph on $\mathcal B$ by joining distinct $P,P'$ if there
is a color set $Q\in\mathcal A$ such that $H$ has an edge between
$Q$ and $P$ and an edge between $Q$ and $P'$. This is a dependency
graph for the events $(E_P)_{P\in\mathcal B}$, and its maximum
degree is at most $(nD)^2$.
\end{lemma}

\begin{proof}
For $P\in\mathcal B$, the event $E_P$ depends only on the
orderings of color sets $Q\in\mathcal A$ that have an edge of $H$
to $P$. If $P'$ is a nonneighbor of $P$ in the proposed dependency
graph, the orderings determining $E_{P'}$ are disjoint from those
determining $E_P$. Since all orderings are independent, $E_P$ is
independent of the joint collection of events associated with its
nonneighbors.

Each color set contains $n$ colors, each of degree at most $D$ in
$H$. Thus at most $nD$ edges of $H$ have an endpoint in any one
color set. In particular, $P$ has edges to at most $nD$ color sets
in $\mathcal A$, and each of those sets has edges to at most $nD$
color sets in $\mathcal B$. Therefore $P$ has at most $(nD)^2$
neighbors in the proposed dependency graph.
\end{proof}

\subsection{Proof of Theorem~\ref{thm:finite}}\label{sec:finite-proof}

\begin{proof}
Recall that each color set has size $n\ge\lceil e^{10}\rceil$. Use the random orderings from Section~\ref{sec:setup}: for each $Q=L(a)$ with $a\in A$, choose an independent uniform bijection $\pi_Q:[n]\to Q$. The color $\pi_Q(i)$ is the color placed in row $i$. Write $\mathcal A=\{L(a):a\in A\}$ and $\mathcal B=\{L(b):b\in B\}$ for the two collections of color sets. Set
\[
 \ell=\log n,\qquad \lambda=\ell-\log\ell-6,\qquad
 D=\Delta(H),\qquad h=\left\lfloor\frac n\ell\right\rfloor.
\]
Then $\ell\ge10$ and the hypothesis is $D\le n\lambda$. The function $t-\log t-6$ is increasing for $t>1$ and is positive at $t=10$, so $\lambda>0$. Also $\lambda<\ell$ because $\log\ell+6>0$.

The proof proceeds as follows. We first give a lower bound on $\tau=(1-1/(n-h+1))^D$, the parameter in Lemma~\ref{lem:rectangle}. We then sum the bounds for all-zero submatrices to bound the probability that one graph $J_P$ has no perfect matching. Finally, we count how many such events can depend on a common ordering and apply Lemma~\ref{lem:local}. The estimates below verify each of these steps.

Since $\ell\ge10$, we have $n/\ell=e^\ell/\ell>2$. Thus
$h=\lfloor n/\ell\rfloor$ satisfies
\begin{equation}\label{eq:h-range}
 \frac{n}{2\ell}\le h\le\frac{n}{\ell}\le\frac n2.
\end{equation}
In particular, $1\le h\le\lfloor n/2\rfloor$.

To estimate $\tau$, use $-\log(1-z)\le z/(1-z)$ with
$z=1/(n-h+1)$. By $D\le n\lambda$ and \eqref{eq:h-range},
\begin{equation}\label{eq:tau}
 \begin{aligned}
 \tau
 &\ge \exp\left(-\frac{D}{n-h}\right)
 \ge \exp\left(-\frac{\lambda\ell}{\ell-1}\right)\\
 &\ge e^{-\lambda-2}
 =\frac{e^4\ell}{n}.
 \end{aligned}
\end{equation}
The last inequality uses $\lambda<\ell$ and $\ell\ge10$.

Fix a color set $P=L(b)$ with $b\in B$. Recall that $J_P$ has sides $[n]$ and $P$, with an edge $ic$ when color $c$ is available in row $i$. Let $E_P$ be the event that $J_P$ has no perfect matching. By Lemma~\ref{lem:hall-obstruction}, this happens exactly when some submatrix $M_P[S,X]$ with $|S|+|X|=n+1$ is all zero. Each such all-zero event is denoted by $\mathcal R(S,X)$. Applying the union bound to the event equality \eqref{eq:hall-failure-events} gives
\[
 \Prob(E_P)\le
 \sum_{\substack{S\subseteq[n],\ X\subseteq P\\ |S|+|X|=n+1}}
 \Prob(\mathcal R(S,X)).
\]
In this sum, both $S$ and $X$ are nonempty, since
$|S|,|X|\le n$ and $|S|+|X|=n+1$. Write $\mu=\min\{|S|,|X|\}$ and $\nu=\max\{|S|,|X|\}$. Since $\mu+\nu=n+1$ and $\nu\ge \mu$, we have $\nu\ge(n+1)/2\ge n/2$. We split the sum according to whether $\mu\le h$ or $\mu>h$, because the number of columns or rows used in Lemma~\ref{lem:rectangle} is $\min\{\mu,h\}$. This number is $\mu$ in the first range and $h$ in the second. In the first range we count submatrices separately for each $\mu$. In the second, a single probability bound will apply to all submatrices, allowing us to use the total count $4^n$.

For $\mu\le h$, we have $\nu\ge n/2$ and, by \eqref{eq:tau},
\[
 \nu\mu\tau\ge\frac{e^4}{2}\mu\ell\ge27\mu\ell.
\]
Hence Lemma~\ref{lem:rectangle} gives
\[
 \Prob(\mathcal R(S,X))\le n^{-27\mu}.
\]
For each fixed $\mu\le h$, Lemma~\ref{lem:submatrix-count} gives at
most $2n^{2\mu-1}$ pairs $(S,X)$. Their total contribution is therefore
at most
\begin{equation}\label{eq:small}
 \sum_{\mu=1}^{h}2n^{2\mu-1}n^{-27\mu}
 \le\frac{2n^{-26}}{1-n^{-25}}
 \le4n^{-26}.
\end{equation}

For $\mu>h$, using \eqref{eq:h-range} and \eqref{eq:tau}, we obtain
\[
 \nu h\tau\ge\frac n2\frac{n}{2\ell}\frac{e^4\ell}{n}
 =\frac{e^4n}{4}.
\]
There are at most $4^n$ pairs $(S,X)$ in total by
Lemma~\ref{lem:submatrix-count}, so the contribution from this range
is at most
\begin{equation}\label{eq:large}
 4^n e^{-e^4n/4}\le e^{-12n},
\end{equation}
where $e^4/4-\log 4>12$.

Combining \eqref{eq:small} and \eqref{eq:large}, we obtain, for every
$P\in\mathcal B$,
\begin{equation}\label{eq:failure}
 \Prob(E_P)\le p:=4n^{-26}+e^{-12n}.
\end{equation}

We apply Lemma~\ref{lem:local} using the probability bound
\eqref{eq:failure}. By Lemma~\ref{lem:matching-dependencies}, the
events $E_P$ have a dependency graph of maximum degree
$d\le(nD)^2$. Since
$D\le n\lambda<n\ell$ and $\ell^2\le e^\ell=n$ for $\ell\ge10$,
\begin{equation}\label{eq:dependency-degree}
 d\le(nD)^2\le n^4\ell^2\le n^5.
\end{equation}

We verify the condition of Lemma~\ref{lem:local} using
\eqref{eq:failure} and \eqref{eq:dependency-degree}. Since
$\log n\le n$, we have $e^{-12n}\le n^{-12}$. Using $e<3$ and
$n^5+1\le2n^5$, we obtain
\begin{equation}\label{eq:local-lemma-condition}
 \begin{aligned}
 ep(d+1)
 &\le 6n^5(4n^{-26}+e^{-12n})\\
 &\le 24n^{-21}+6n^{-7}<1,
 \end{aligned}
\end{equation}
where the last inequality holds for $n\ge22027$.

By \eqref{eq:local-lemma-condition}, Lemma~\ref{lem:local} applies. It shows that some choice of orderings makes every graph $J_P$ have a perfect matching. Fix these orderings. By \eqref{eq:packing-matchings}, the rows can then be completed to a full packing. Each row is an independent transversal, and every color occurs in exactly one row, so this gives the partition into $n$ independent transversals required by Theorem~\ref{thm:finite}.
\end{proof}

\subsection{Proof of Theorem~\ref{thm:main}}\label{sec:main-proof}

We now choose the size of each color set as a function of the maximum degree of $H$ and apply Theorem~\ref{thm:finite}.

\begin{proof}
Let $D\ge\lceil e^{20}\rceil$, and let $H$ be a correspondence
cover satisfying the hypotheses of Theorem~\ref{thm:main}. Put
\[
 T=\log D,\qquad
 \beta=T-2\log T-6,\qquad
 n=k(D)=\left\lceil\frac{D}{\beta}\right\rceil.
\]
Since $T\ge20$, we have $2<\beta\le T$. The function $e^t/t$ is
increasing for $t\ge20$, so
\begin{equation}\label{eq:main-size-check}
 n\ge\frac{D}{\beta}\ge\frac{D}{T}
 \ge\frac{e^{20}}{20}>e^{10}.
\end{equation}
Thus $n\ge\lceil e^{10}\rceil$. Also,
\[
 n\le\frac{D}{\beta}+1\le\frac D2+1\le D.
\]

The bounds $D/T\le n\le D$ give
\[
 \log n\ge T-\log T,
 \qquad
 \log\log n\le\log T.
\]
Hence
\[
 \log n-\log\log n-6\ge T-2\log T-6=\beta.
\]
Since $n\ge D/\beta$ and $\Delta(H)\le D$, we obtain
\begin{equation}\label{eq:main-degree-check}
 n(\log n-\log\log n-6)\ge n\beta\ge D\ge\Delta(H).
\end{equation}
Theorem~\ref{thm:finite} therefore gives a partition into
$n=k(D)$ independent transversals.

Finally, write $n=D/\beta+\theta_D$, where $0\le\theta_D<1$. Then
\[
 \frac{k(D)}{D/\log D}
 =\frac{T}{T-2\log T-6}
  +\theta_D\frac{T}{D}
 \longrightarrow1
 \qquad(D\to\infty).
\]
We are done.
\end{proof}

\section*{Conflict of interest}
The author declares no conflict of interest.

\section*{Data availability}
No data were used in this study.

\section*{Declaration of AI use}
The author used GPT-6.1 Sol (OpenAI) to assist with language editing
and improve the clarity of the manuscript. The author independetly wrote the
initial version of the paper using $k(D)=(1+\varepsilon)D/\log D$
throughout the proof. After discussing that version with the author,
GPT-6.1 Sol suggested the explicit choice
\[
 k(D)=\left\lceil
 \frac{D}{\log D-2\log\log D-6}
 \right\rceil
\]
and helped the author update the numerical estimates throughout the
paper accordingly. The author reviewed the resulting text, rewrote
it in his own words, and takes full responsibility for the content.

\end{document}